\documentclass{amsart}

\newtheorem{teo}{Theorem}[section]
\newtheorem{prop}[teo]{Proposition}
\newtheorem{cor}[teo]{Corollary}
\newtheorem{lema}[teo]{Lemma}
\newtheorem{claim}[teo]{Claim}

\theoremstyle{definition}
\newtheorem{dfn}[teo]{Definition}
\newtheorem{ex}[teo]{Example}

\theoremstyle{remark}
\newtheorem{rem}[teo]{Remark}

\numberwithin{equation}{section}

\newcommand{\WSt}{\ensuremath{W_{\widetilde{\Sigma}}}}

\newcommand{\WSrt}{\ensuremath{W_{\widetilde{\Sigma}_{r}}}}

\newcommand{\s}{\ensuremath{\varsigma}}

\newcommand{\tub}{\ensuremath{\mathrm{Tub} }}

\newcommand{\F}{\ensuremath{\mathcal{F}}}

\newcommand{\TO}{\ensuremath{\mathcal{O}}}
\newcommand{\singularF}{\ensuremath{\mathcal{X}_{\mathcal{F}}}}

\newcommand{\dank}{\textsf{Acknowledgments.\ }}

\newcommand{\Iso}{\ensuremath{\mathrm{Iso}}}

\usepackage{amssymb}
\usepackage{amsmath}

\newcommand{\R}{\mbox{$\mathbf R$}}

\begin{document}

\title[Singular holonomy of S.R.F.S]{ Singular holonomy of singular riemannian foliations with sections}

\author{Marcos M. Alexandrino} 
 
\address{Marcos Martins Alexandrino \\
Instituto de Matem\'{a}tica e Estat\'{\i}stica,
Universidade de S\~{a}o Paulo (USP),
Rua do Mat\~{a}o, 1010, 
S\~{a}o Paulo, SP 05508-090, Brazil}

\email{malex@ime.usp.br}

\email{marcosmalex@yahoo.de}


\date{April 2007}

\thanks{The author was  supported by CNPq and FAPESP.}


\begin{abstract}
In this paper we review some author's results about singular holonomy of singular riemannian foliation with sections (s.r.f.s for short) and also some results of a joint work with T\"{o}ben and  a joint work with Gorodski. We stress here that the condition that the leaves are compact, used in some of these results, can be replaced by the condition that the leaves are closed embedded.
 We  also briefly recall some of T\"{o}ben's results about blow-up of s.r.f.s. Then
we  use this technique to get conditions under which a holonomy map can be extended to a global isometry.

\end{abstract}


\maketitle

\section{Introduction}
 
A singular riemannian foliation $\F$ on a complete riemannian 
manifold $M$ is said to admit sections if 
each regular point of $M$ is contained in a 
complete totally geodesic immersed 
submanifold $\Sigma$ that meets every leaf of $\F$ orthogonally 
and whose 
dimension is the codimension of the regular leaves of $\F$ (see Definition \ref{dfn-srfs}).

Typical examples of a singular riemannian foliation with section (s.r.f.s for short) are the partition formed by the orbits of a polar action,  
partition formed by parallel submanifolds  of an isoparametric submanifold and partition formed by parallel submanifolds  of an equifocal submanifold
(see definitions in Section 2). Others examples can be constructed by suspension of homomorphism, suitable changes of metric and surgery. 

The property that s.r.f.s are equifocal (see Definition \ref{dfn-equifocal} and Theorem \ref{frss-eh-equifocal}) allow us to extend the (regular) normal holonomy map of regular leaves to the so called \emph{singular holonomy map} which is also defined on  singular points (see Proposition \ref{prop-holonomia-singular}).

Singular holonomy and the Slice Theorem (see Theorem \ref{sliceteorema}) give us a complete description of a s.r.f.s on a neighborhood of a singular point. 

On the other hand, singular holonomy turned out to be a useful tool to study global properties of s.r.f.s (e.g. see Corollary \ref{cor-estrato-trivial-folhas-fechadas}, Theorem \ref{teo-fundamental-domain} and Theorem \ref{teo-s.r.f.s-transnormal}).

In this paper we review some author's results about singular holonomy and also some results of a joint work with T\"{o}ben and  a joint work with Gorodski. We stress here that the condition that the leaves are compact, used in some of these results, can be replaced by the condition that the leaves are closed embedded (see Theorem \ref{teo-holonomia-trivial}, Theorem \ref{teo-fundamental-domain} and Theorem \ref{teo-s.r.f.s-transnormal} ). 
From this new remark we also infer a new result (see Corollary \ref{cor-holonomia-grupo-fundamental}).
Finally, we  briefly recall some of T\"{o}ben's results about blow-up of s.r.f.s (see Theorem \ref{Th-Toeben} and Proposition \ref{Prop-Toeben}). Then
we  use this technique to get conditions under which a holonomy map can be extended to a global isometry (see Proposition \ref{Prop-extensao}).

This paper is organized as follows. In Section 2, we review some facts about s.r.f.s and fix the notation. In Section 3 we recall an author's result conerning to Molino's conjecture and orbits of Weyl pseudogroups (see Theorem \ref{teo-Molino-conjecture}). 
In Section 4 we discuss the relation between the holonomy of a s.r.s.f and the fundamental group of the space. In particular, we prove  Theorem \ref{teo-holonomia-trivial} and reformulate previous results. Finally,  
in Section 5 we review  T\"{o}ben's results about blow-up of s.r.f.s and prove Proposition \ref{Prop-extensao}.

 
\dank The author is grateful to Dirk T\"{o}ben  
for useful suggestions. 
The author acknowlegde the support of  CNPq (Conselho Nacional de Desenvolvimento Cient\'{\i}fico  e Tecnol\'{o}gico - Brasil) and FAPESP(Fundan\c{c}\~{a}o de Amparo \`{a} Pesquisa do Estado de S\~{a}o).


\section{Facts about s.r.f.s.}

In this section, we recall some results about s.r.f.s.\ that will be used 
in this text. Details can be found in~\cite{Alex2}. 
Throughout this section, we assume that $\F$ is a 
singular riemannian foliation with sections
on a complete riemannian manifold $M$; 
we start by recalling its definition.

\begin{dfn}
\label{dfn-srfs}
A partition $\F$ of a complete riemannian manifold $M$ by connected 
immersed submanifolds (the \emph{leaves}) is called a \emph{singular 
riemannian foliation with sections} of $M$ (\emph{s.r.f.s.}, for short) if it 
satisfies
the following conditions:
\begin{enumerate}
\item $\F$ is \emph{singular},
i.e.~the module $\singularF$ of smooth vector fields on $M$ that are 
tangent at each point to the corresponding leaf acts transitively on each 
leaf. In other words, for each leaf $L$ and each $v\in TL$ with 
footpoint $p,$ there exists $X\in \singularF$ with $X(p)=v$.
\item  The partition is \emph{transnormal}, i.e.~every 
geodesic that is perpendicular to a leaf at one point remains 
perpendicular to every leaf it meets.
\item For each regular point $p$, the set $\Sigma :=\exp_{p}(\nu_p 
L_{p})$ 
is a complete immersed submanifold that meets all the leaves and meets
them always orthogonally. The set  
$\Sigma$ is called a \emph{section}.
\end{enumerate}
\end{dfn}

\begin{rem}
The concept of s.r.f.s was introduced in \cite{Alex2} and continued to be studied by me in \cite{Alex1,Alex3,Alex4}, by T\"{o}ben in \cite{Toeben,Toeben2},  by T\"{o}ben and I in \cite{AlexToeben}, by Lytchak and Thorbergsson in \cite{LytchakThorbergsson} and recently by Gorodski and I in \cite{AlexGorodski}. 
In \cite{Boualem} Boualem dealt with a 
singular riemannian foliation $\F$ on a complete manifold $M$ 
such that the distribution of normal spaces of the regular leaves  
is integrable. It was proved in \cite{Alex4} that such an 
$\F$ must be a s.r.f.s.\/ and, in addition,
the set of regular points is open and dense in each section. 
\end{rem}

A typical example of s.r.f.s is the partition formed by the orbits of a polar action. An isometric action of a compact Lie group $G$ on a complete riemannian manifold $M$ is called \emph{polar} if there exists a complete immersed submanifold $\Sigma$ of $M$ that meets all $G$-orbits orthogonally and whose dimension is equal to the codimension of a regular orbit.

Another typical example of a s.r.f.s is the partition formed by parallel submanifolds  of an isoparametric submanifold $N$ of an euclidean space.
A submanifold $N$ of an euclidean space is called \emph{isoparametric} 
if its normal bundle is flat and the principal curvatures along any 
parallel 
normal vector field are constant. 
Theorem~\ref{sliceteorema} below shows how s.r.f.s.\ and 
isoparametric foliations are related to each other.
In order to state this theorem, we need the concepts of slice and  
local section. 
Let $q\in M$, and let $\tub(P_{q})$ be a tubular neighborhood of a 
plaque $P_{q}$ that contains $q$. Then the connected component of 
$\exp_{q}(\nu P_{q})\cap \tub(P_{q})$  that contains $q$ is called a 
\emph{slice} at $q$ and is usually denoted by $S_{p}.$ 
A \emph{local section} $\sigma$ (centered at $q$) of a section $\Sigma$  is a connected component $\tub(P_{q})\cap\Sigma$ (which contains $q$).  
 
\begin{teo}[\cite{Alex2}]
\label{sliceteorema}
Let $\F$ be a s.r.f.s. on a complete riemannian manifold $M.$ 
Let $q$ be a singular point of $M$ and let 
$S_{q}$ a slice at $q.$ Then
\begin{enumerate}
\item  Denote  
$\Lambda(q)$  the set of 
local sections $\sigma$ centered at $q$ 
 Then  $S_{q}= \cup_{\sigma\in\Lambda (q)}\, \sigma$.
\item $S_{x}\subset S_{q}$ for all $x\in S_{q}$.
\item $\F|S_q$ is a s.r.f.s.\ on $S_{q}$ with the induced metric from 
$M$.
\item $\F|S_q$ is diffeomorphic to an isoparametric foliation on an 
open 
subset of $\mathbf{R}^{n}$, where $n$ is the dimension of $S_{q}$.
\end{enumerate}
\end{teo}

From (d), it is not difficult to derive the following corollary.

\begin{cor}
\label{estratificacao-singular}
Let $\sigma$ be a local section. Then the set of singular points of 
$\F$ 
that are contained in $\sigma$ is a finite union of totally 
geodesic hypersurfaces. 
These hypersurfaces are mapped by a diffeomorphism to the focal 
hyperplanes 
contained in a section of an isoparametric foliation on an open subset 
of an euclidean space.   
\end{cor}

We will call the set of singular points of $\F$ contained in $\sigma$ 
the \emph{singular stratification of the local section} $\sigma$. 
Let $M_{r}$ denote the set of regular points in $M.$ A \emph{Weyl 
Chamber} of 
a local section $\sigma$ is the closure in $\sigma$ of a connected 
component 
of $M_{r}\cap\sigma$. One can prove that a Weyl Chamber of a 
local section is a convex set.  

 Theorem~\ref{sliceteorema} also implies that a s.r.f.s can be locally trivialized by a transnormal map, whose definition we recall now.  

\begin{dfn}[Transnormal  Map]
\label{dfn-transnormal-map}
Let  $M^{n+q}$  be a complete riemannian manifold. A smooth map   
$F=(f_{1},\ldots, f_{q}):M^{n+q}\rightarrow \mathbf{R}^{q}$ is 
called a \emph{transnormal map} if the following assertions hold:
\begin{enumerate}
\item[(0)]  $F$ has a regular value.
\item[(1)]  For each regular value $c$ there exists a neighborhood  $V$ of $ F^{-1}(c)$ in $M$ such that $F \mid_{V}\rightarrow F(V)$ is an integrable riemannian submersion.
\end{enumerate}
In particular, a transnormal map $F$ is said to be an \emph{isoparametric map} if  $V$ can be chosen to be $M$  and   $\bigtriangleup f_{i}= a_{i}\circ F,$ where $a_{i}$ are smooth functions. 
\end{dfn}

\begin{rem}
Recall that each isoparametric submanifold in an euclidian space
can always be described as a regular level set of an isoparametric  
polynomial map (see Terng\cite{terng} or Terng and Palais \cite{PTlivro}). 
On the other hand, the regular leaves 
of an analytic transnormal map on a complete analytic
riemannian manifold are equifocal manifolds and leaves of a 
s.r.f.s (see ~\cite{Alex1}).
\end{rem}

\begin{prop}
\label{prop-placas-transnormais}
The plaques of a s.r.f.s. are always level sets of a transnormal map. 
\end{prop}

In \cite{TTh1}, Terng and Thorbergsson introduced the concept of 
equifocal submanifolds with flat sections in symmetric spaces in 
order to generalize the definition of isoparametric submanifolds 
in euclidean space. 
Next we review the slightly more general
definition of equifocal submanifolds in riemannian 
manifolds.

\begin{dfn}
\label{dfn-equifocal}
A connected immersed submanifold $L$ of a complete riemannian manifold 
$M$ is called \emph{equifocal} if it satisfies the following 
conditions:
\begin{enumerate}
\item The normal bundle $\nu(L)$ is flat.
\item $L$ has sections, i.e.~for each~$p\in L$, 
the set $\Sigma :=\exp_{p}(\nu_p L_{p})$ is a complete immersed totally 
geodesic submanifold.
\item For each parallel normal field $\xi$ on a neighborhood $U \subset L$,  
the derivative of the map $\eta_{\xi}:U\to M$ defined 
by $\eta_{\xi}(x):=\exp_{x}(\xi)$ has constant rank.
\end{enumerate}
\end{dfn}

The next theorem relates s.r.f.s.\ and equifocal submanifolds.

\begin{teo}[\cite{Alex2}]
\label{frss-eh-equifocal}
Let $L$ be a regular leaf of a s.r.f.s.\ $\F$ of a 
complete riemannian manifold $M$.
\begin{enumerate}
\item Then $L$ is equifocal. In particular, the union of the 
regular leaves that 
have trivial normal holonomy is an open and dense set in $M$ provided 
that all the leaves are compact.
\item Let $\beta$ be a smooth curve of $L$ and $\xi$ a parallel normal 
field to
$L$ along~$\beta$. 
Then the curve $\eta_{\xi}\circ \beta$ belongs to a leaf of $\F$.
\item Suppose that $L$ has trivial holonomy and let $\Xi$ denote the 
set of 
all parallel normal fields on $L$. 
Then $\F=\{\eta_{\xi}(L)\}_{\xi\in \, \Xi}.$ 
\end{enumerate}
\end{teo}

The above theorem allows us to define the 
singular holonomy map, which will be very useful to study $\F$.

\begin{prop}[Singular holonomy map] 
\label{prop-holonomia-singular}
Let $\F$ be a s.r.f.s. on a complete riemannian manifold $M$ and $q_0$ and $q_1$ two points contained in a leaf $L_{q}.$ Let  $\beta:[0,1]\rightarrow L_{p}$ be a smooth curve contained in a regular leaf $L_{p}$, such that $\beta(i)\in S_{q_{i}},$ where $S_{q_{i}}$ is the slice at $q_{i}$ for $i=0,1.$ Let $\sigma_{i}$ be a local section contained in $S_{q_{i}}$ which contains $\beta(i)$ and $q_{i}$ for $i=0,1$.  Finally let $[\beta]$ denote the homotopy class of $\beta.$ Then there exists an isometry $\varphi_{[\beta]}:U_{0} \rightarrow U_{1},$ where the source $U_{0}$ and target $U_{1}$ are contained in $\sigma_{0}$ and $\sigma_{1}$ respectively,  which has the following properties:
\begin{enumerate}
\item[1)] $q_{0}\in U_{0}$
\item[2)]$\varphi_{[\beta]}(x)\in L_{x}$ for each $x\in U_{0}.$
\item[3)]$d\varphi_{[\beta]}\xi(0)=\xi(1),$ where $\xi(s)$ is a  parallel normal field along $\beta(s).$
\end{enumerate} 

\end{prop}

An isometry as in the above proposition
is called the \emph{singular holonomy map along $\beta$}. 
 
We remark that, in the definition of the singular holonomy map, 
singular points can be contained in the domain $U_0.$  
If the domain $U_0$ and the range $U_1$ are sufficiently small, then the 
singular holonomy map coincides with the usual holonomy map along 
$\beta$.

Theorem \ref{sliceteorema} establishes a relation between s.r.f.s.\ and 
isoparametric foliations. Similarly as in 
the usual theory of isoparametric submanifolds, it is 
natural to ask if we can define a (generalized) Weyl group action on 
$\sigma$. The following definitions and results 
deal with this question.

\begin{dfn}[Weyl pseudogroup $W$]
\label{definitionWeylPseudogroup}
 The pseudosubgroup generated by all singular holonomy maps 
$\varphi_{[\beta]}$ such that $\beta(0)$ and $\beta(1)$ belong to the 
same 
local section~$\sigma$ is called the \emph{generalized Weyl 
pseudogroup}
of~$\sigma$. Let $W_{\sigma}$ denote this pseudogroup. 
In a similar way, we define $W_{\Sigma}$ for a section $\Sigma$. 
Given a slice~$S$, we define~$W_{S}$ as the set of all singular 
holonomy 
maps~$\varphi_{[\beta]}$ such that~$\beta$ is contained in the 
slice~$S$.
\end{dfn}

\begin{rem}
Regarding the definition of pseudogroups and orbifolds,
see Salem~\cite[Appendix D]{Molino}. 
\end{rem}

\begin{prop}
\label{propWisinvariant}
Let $\sigma$ be a local section. Then the reflections in the 
hypersurfaces of the singular stratification of the local section $\sigma$ leave 
$\F|\sigma$  invariant. Moreover these reflections are elements of 
$W_{\sigma}.$
\end{prop}

By using the technique of suspension, 
one can construct an example of a s.r.f.s.\ 
such that  $W_{\sigma}$ is larger than the pseudogroup 
generated by the reflections in the hypersurfaces of the 
singular stratification of $\sigma$. On the other hand,
a sufficient condition 
to ensure that both pseudogroups coincide is that the leaves of $\F$ 
have trivial normal holonomy and be compact.
So it is natural to ask under which conditions we can garantee that the normal holonomy of regular leaves are trivial. This question will be answered in Section 4.


\section{Molino's conjecture and singular holonomy}

In this section we review an author's result concerning to Molino's conjecture. 

In \cite{Molino} Molino proved that, if $M$ is compact, the closure of the leaves of a (regular) riemannian foliation form a partition of $M$ which is a singular riemannian foliation. He also proved that the leaf closure are orbits of a locally constant sheaf of germs of (transversal) Killing fields. 
If the foliation is a singular riemannian foliation and $M$ is compact, then Molino was able to prove (see \cite{Molino} Theorem 6.2 page 214) that the closure of the leaves should be a transnormal system, but as he remarked, it remains to prove that the closure of the leaves is in fact a singular foliation. In \cite{Alex4} I proved the Molino's conjecture, when $\F$ is a s.r.f.s. In addition I studied the singular holonomy of $\F$ and in particular the   tranverse orbits of the closure of a leaf. In this work was not assumed that  $M$ should be compact.

\begin{teo}[\cite{Alex4}]
\label{teo-Molino-conjecture}
Let $\F$ be a s.r.f.s. on a complete riemannian manifold $M.$ 
Then
\begin{enumerate}
\item[a)] the closure of the leaves of $\F$ form a partition of $M$ which is a singular riemannian foliation,i.e, $\{\overline{L}\}_{L\in\F}$ is a singular riemannian foliation.
\item[b)] Each point $q$ is contained in an homogenous submanifold $\TO_{q}$ (possible with dimension $0).$ If we fix a local section $\sigma$ that contains $q,$ then $\TO_{q}$ is a connected component of an orbit of the closure of the Weyl pseudogroup of $\sigma.$   
\item[c)] If $q$ is a point of the submanifold $\overline{L},$ then a neighborhood of $q$ in $\overline{L}$ is the product of the homogenous submanifold  $\TO_{q}$ with plaques with the same dimension of the plaque $P_{q}.$
\item[d)] Let $q$ be a singular point and  $T$  the intersection of the slice $S_{q}$ with the   singular stratum that contains $q.$ Then the normal connection of $T$ in $S_{q}$ is flat.
\item[e)] Let $q$ be a singular point and $T$ defined as in Item d). Let $v$ be a parallel normal vector field along $T,$ $x\in T$ and $y=\exp_{x}(v)$. Then  
$\TO_{y}=\eta_{v}(\TO_{x}).$ 
\end{enumerate}
\end{teo}

One can construct examples that illustrate the above theorem by means of suspension of homomorphisms (see Example \ref{ex-molino-conjecture}). In fact, the suspension technique is very useful to construct examples of s.r.f.s. with nonembedded leaves, with exceptional leaves and also inhomogeneous examples. Other techniques to construct examples of s.r.f.s on nonsymmetric spaces are  suitable changes of metric and  surgery (see \cite{AlexToeben} for details). 

\begin{cor}
\label{cor-estrato-trivial-folhas-fechadas}
Let $\F$ be a s.r.f.s. on a complete manifold $M$ and $q$ a singular point. Let $T$ denote the intersection of the slice $S_{q}$ with the  stratum that contains $q.$ Suppose that $T=\{q\}.$ Then all the leaves of $\F$ are closed.
\end{cor}

\begin{rem}
According to the slice theorem (see Theorem \ref{sliceteorema})  the restriction of the foliation $\F$ to the slice $S_{q}$ is diffeomorphic to an isoparametric foliation $\widetilde{\F}$ on an open set of an euclidean space. Therefore the condition that $T$ is a point is equivalent to saying that a regular leaf of $\widetilde{\F}$ is a full isoparametric submanifold.
\end{rem}

\begin{ex}
\label{ex-molino-conjecture}
In what follows we construct a s.r.f.s  such that the intersection of a local section with the closure of a regular leaf is  an orbit of an action of a subgroup of isometries of the local section. This isometric action is not a polar action. This implies that \emph{there exists a s.r.f.s $\F$ such that the partition formed by the closure of the leaves of $\F$ is a singular riemannian foliation without sections}.
\end{ex}

 First, define an homomorphism $\rho$ as 

\[
\begin{array}{llll}
\rho:& \pi_{1}(\mathbf{S}^{1},q_{0})&\rightarrow& \Iso (\mathbf{R}^{2}\times\mathbf{C}\times\mathbf{C})\\
     &  n               &\rightarrow&((x,z_{1},z_{2})\rightarrow (x, e^{(i\,n\,k)}\cdot z_{1}, e^{(i\,n\,k)}\cdot z_{2} ))\\
    \end{array}\]
where $k$ is an irrational number.
Set $\widetilde{M}:=\mathbf{R}\times (\mathbf{R}^{2}\times\mathbf{C}\times\mathbf{C})$
and define an action of $\pi_{1}(\mathbf{S}^{1},b_{0})$ on $\widetilde{M}$ by
\[ [\alpha]\cdot(\hat{b},t):=([\alpha]\cdot\hat{b},\rho(\alpha^{-1})\cdot t),\]
where $[\alpha]\cdot \hat{b}$ denotes  the  deck transformation associated to $[\alpha]$ applied to a point $\hat{b}\in\mathbf{R}.$ 
Let $M:=\widetilde{M}/\sim$ be the orbit space, 
  $\Pi:\widetilde{M}\rightarrow M$ the canonical projection and 
 $\hat{P}:\mathbf{R}\rightarrow \mathbf{S}^{1}$ the covering map.
Finally define a map $P:M\rightarrow \mathbf{S}^{1}$ as follows

\[ \begin{array}{rlcl}
P: &M &\rightarrow& \mathbf{S}^{1}\\
   &\Pi(\hat{b},t)& \rightarrow & \hat{P}(\hat{b})
\end{array}
\]

It is possible to prove that  
 $M$ is a total space of a fiber bundle, which has $P$ as the   projection over the basis $\mathbf{S}^{1}.$ Besides the fiber of this bundle is  $\mathbf{R}^{2}\times\mathbf{C}\times\mathbf{C}$ and the structural group is given by the image of $\rho.$

Now let  $\hat{\F}_{0}$ be the singular foliation of codimension 5 on 
$\mathbf{R}^{2}\times\mathbf{C}\times\mathbf{C}$
 whose leaves are the product of  points in $\mathbf{C}\times\mathbf{C}$ with  circles  in $\mathbf{R}^{2}$ centered at $(0,0).$  It is easy to see that the foliation $\hat{\F}_{0}$ is a singular riemannian foliation with sections.

Finally  set $\F:=\Pi(\mathbf{R}\times \hat{\F_{0}}).$ It turns out  that $\F$ is s.r.f.s. such that the intersection of the section $\Pi(0\times\mathbf{R}\times\mathbf{C}\times\mathbf{C})$
with the closure of a regular leaf is an orbit of an isometric action on the section. 
This isometric action is not a polar action, since the isometric action 
\[
\begin{array}{lcl}
 \mathbf{S}^{1}\times\mathbf{C}\times\mathbf{C}&\rightarrow & \mathbf{C}\times\mathbf{C}\\
    (s,z_{1},z_{2}) &\rightarrow& (s\cdot z_{1},s\cdot z_{2})\\
    \end{array}\]
is not a polar action.


\section{holonomy and fundamental group of $M$}

In this section we review a result of a joint work with T\"{o}ben \cite{AlexToeben}, where we discuss the relation between holonomy of a s.r.f.s $\F$ on $M$ and the fundamental group of $M$. We stress here that the condition that the leaves are compact, which was used in \cite{AlexToeben}, can be replaced by the condition that the leaves are embedded and closed (see Theorem \ref{teo-holonomia-trivial}). This remark allow us to generalize a previus  result (see Theorem \ref{teo-fundamental-domain}) and  infer a new one(see Corollary \ref{cor-holonomia-grupo-fundamental}). We also briefly recall a result of a joint work with Gorodski \cite{AlexGorodski} (see Theorem \ref{teo-s.r.f.s-transnormal}).



\subsection{Transversal Frame Bundle of a s.r.f.s}
\label{subsection-transversal-frame-bundle}

In \cite{Molino} Molino associated an $O(k)$-principal bundle, the {\it orthogonal transverse frame bundle} to a regular riemannian foliation $(M,\F)$ of codimension $k$. A fiber of this bundle over a point $p$ in $M$ is defined as the set of orthonormal $k$-frames in $\nu_p L_p$, where $L_p$ is the leaf through $p$. Proposition \ref{prop-fibrado-referenciais} generalize this notion for a s.r.f.s. of codimension $k$. Its restriction to the regular stratum $M_r$ will coincide with the orthogonal transverse frame bundle in the sense of Molino.  We will use this bundle to review the proof of  Theorem \ref{teo-holonomia-trivial} and Theorem \ref{teo-holonomia-trivial-secao-flat}.

\begin{prop}[\cite{AlexToeben}]
\label{prop-fibrado-referenciais}
Let $\F$ be a s.r.f.s on a complete riemannian manifold $M.$ 
\begin{enumerate}
\item[a)] There exists a continuous principal $O(k)$-bundle $E$ over $M$ that is associated to the s.r.f.s. $\F.$ The restriction of $E$ over $M_{r}$ (denoted by $E_{r}$) coincides with the usual orthogonal transverse frame bundle of the riemannian foliation $\F_{r}$(the restriction of $\F$ to $M_{r}$).
\item[b)] There exists a singular $C^0$-foliation $\widetilde{\F}$ on $E.$ The restriction of $\widetilde{\F}$ to $E_{r}$ coincides with the usual parallelizable foliation $\widetilde{\F}_{r}$ on $E_{r},$ which is a foliation with trivial holonomy whose leaves cover the leaves of $\F_{r}.$
\item[c)] There exist $C^{0}$ holonomy map associated to $\widetilde{\F}$; hence we can define a Weyl pseudogroup $\WSt$.
\item[d)] There exist $C^{0}$local trivializations of $\widetilde{\F}$.
\item[e)] If the sections of $\F$ are flat, the bundle $E$, the foliation $\widetilde{\F}$, holonomy maps and trivializations are smooth.
\end{enumerate}
\end{prop}     

Before we sketch the construction of the bundle $E$, we present a very simple example.

\begin{ex}
\label{trivialex} 
Consider  $M:=\mathbf{R}^{2}$ foliated by circles centered at the origin. We denote it by $\F$. The only singular leaf is the origin and the sections are the lines through the origin. Excising the singular leaf we obtain a regular riemannian foliation $\F_r$ of $M_{r}:=\mathbf{R}^{2}-\{(0,0)\}$. 
Let $E_{r}$ be the orthogonal transverse frame bundle (in the sense of Molino) associated to $\F_{r}.$ It is not difficult to see that $E_{r}=M^{1}_{r} \amalg M^{-1}_{r},$ where $M^{i}_{r}:=(\mathbf{R}^{2}-\{(0,0)\})\times \{i\}$  for $i=1,-1.$ We can identify $M^{1}_{r}$ (respectively $M^{-1}_{r}$) with the unit normal field outward (respectively inward) oriented.
 Set  $E= M^{1}\cup M^{-1},$ where $M^{i}:= M_{r}^{i}\cup (\{(0,0)\}\times \{i\}).$ We  will define $E$ as the  transverse frame bundle associated to $\F.$ It is obvious that the restriction of $E$ to $\pi^{-1}(M_r)$ is the orthogonal tranverse frame bundle $E_{r}.$ The set $E$ can also be regarded as a set of equivalence classes, where the equivalence is defined as follows. 
Let $(\zeta^{i}_{p},C_{i})$ for $i=1,2$ be a pair of a vector  tangential to some local sections $\sigma$ through $p$ with footpoint $p\in M$ and Weyl Chamber $C_{i}$  in $\sigma$ that contains $p.$ Then $(\zeta^{i}_{p},C_{i})$ are defined to be equivalent if there exists a rotation $\varphi$ (a holonomy map) such that $\varphi (C_{1})=\varphi (C_{2})$ and $\varphi_*\zeta^{1}=\varphi_*\zeta^{2}.$
We say that $\tilde{p}:=[(\zeta_{p},C)]$ belongs to $M^{1}$ (respectively $M^{-1}$)  if a representative $(\zeta_{p},C)$ induces the outward  (respectively inward) orientation of $\F_{r}$
 by parallel transport along the Weyl chamber and the circles.
 Note that if $p$ is not $(0,0)$ then there exists only one Weyl chamber $C$ that contains $p$ and hence this new definition coincides with the definition of a vector $\zeta_{p}$ with footpoint $p,$ when $p$ is regular.  

\end{ex}

\textsf{Sketch of construction of the Transversal Frame Bundle $E.$  }

Let $\F$ be a s.r.f.s. on a complete riemannian manifold $M.$ 
Let $(\zeta_p,C)$ be a pair of an orthonormal $k$-frame $\zeta$ with footpoint $p$ tangential to a local section $\sigma$ and the germ of a Weyl chamber $C$ of $\sigma$ at $p$.
We identify $(\zeta^1_p,C_1)$ and $(\zeta^2_p,C_2)$ if there is a holonomy map $\varphi\in W_{S_p}$ (which fixes $p$) that maps  $ C_1$ to $C_2$ as germs in $p$ and $\zeta^1_p$ to $\zeta^2_p$ at first order, where  $W_{S_p}$ is the set of all holonomy map $\varphi_{[\beta]}$ such that $\beta$ is contained in the slice  $S_p$. In other words, the equivalence class $[(\zeta_p,C)]$ consists of the $W_{S_p}$-orbit $(\varphi_*\zeta_p,\varphi(C)), \varphi\in W_{S_p}$. We call an equivalence class $[(\zeta_p,C)]$ {\it transverse frame}, and the set $E$ of transverse frames {\it transverse frame bundle}.

Let $\pi:E\to M$ be the footpoint map. The fiber $F_q=\pi^{-1}(q)$ is equal to the set of transverse frames $[(\zeta_q,C)]$. There is a natural right action of $O(k)$ on $E$ by $[(\zeta_q,C)]\cdot g:=[(\zeta_q\cdot g,C)]$. This action is well-defined and simply transitive on the fiber. Note that in each equivalence class there is only one representative with a given  Weyl chamber.

Given a transverse frame $\tilde{q}=[(\zeta_{q},C)]$ over a point $q$ it is possible to use parallel transport and  equifocality of $\F$ to find a neighborhood $U$ of $q$ in $M$ and a map  
$\s:U\rightarrow \pi^{-1}(U)$ such that $\s(q)=\tilde{q}$ and $\pi\circ\s(x)=x$ for $x\in U$. With the cross section $\s$ we can define a  trivialization of $E|U$ as follows. 
\begin{eqnarray}
\label{trivializationphi}
\phi:U\times O(k)&\to& E|U\nonumber\\
 (x,g)&\mapsto & \s(x)\cdot g.
\end{eqnarray}

Let $E|U$ take the induced topology via $\phi$. One has to show that this topology on $E$ is coherently defined, i.e., the transition from one trivialization to another trivialization is a homeomorphism. This follows from the next lemma.
\begin{lema}
\label{lema-funcao-h}
Consider two cross-sections $\s_i:U_i\to E$ with $U_i\cap U_j\neq\emptyset$ and the corresponding trivializations $\phi_i:U_i\times O(k)\to E|U_i$.  Define $h:U_1\cap U_2\to O(k)$ by $\s_1(x)=\s_2(x)\cdot h(x)$. Then 
\begin{enumerate}
\item $\phi_2^{-1}\circ\phi_1(x,g)=(x,h(x)\cdot g)$.
\item $h$ is constant along the plaques in $U_1\cap U_2$.
\item The map $h:U_1\cap U_2\to O(k)$ is continuous at all points and differentiable at all regular points. If the sections are flat, $h$ is locally constant.
\end{enumerate}
\end{lema}
Finally we define a singular foliation $\widetilde{\F}$ on $E$ as follows: Let $\phi:U\times O(k)\to E|U$ be a trivialization and $P_x$ for $x\in U$ the plaque of $\F$ in $U$. We define $\widetilde{\F}|U$ by the partition $\widetilde{P}_{\phi(x,g)}:=\phi(P_x,g)$. Since the transition map $h$ is constant along the plaques, $\widetilde{\F}$ is well-defined on $E$. We define a leaf $\widetilde{L}$ through a point $x$ as the set of endpoints of continuous paths contained in plaques that start in $x$. The restriction of $\widetilde{\F}$ to the bundle $E_r=E|M_r$ over the regular stratum $M_r$ is the standard foliation described in \cite{Molino}.
 

\subsection{Main result}

In \cite{AlexToeben} T\"{o}ben and I proved that the holonomy of the leaves of a s.r.f.s $\F$ on a complete manifold $M$ is trivial, if the leaves are compact and $M$ is simply connected. I stress here that the condition that the leaves are compact is too strong. It is enought to assume that the leaves are embedded and closed. 

\begin{teo}
\label{teo-holonomia-trivial}
Let $\F$ be a s.r.f.s.\ on a simply connected riemannian manifold $M$. 
Suppose also that the leaves of $\F$ are embedded and closed. Then
 each regular leaf has trivial holonomy.
\end{teo}
\begin{proof}

In this proof we need the concept of fundamental group of a pseudogroup, which we briefly recall below (for details see Salem  \cite[Appendix D]{Molino}).

We start by recalling the definition of a $W$-loop of a pseudogroup $W$ on a $C^{0}$ manifold $\Sigma.$ 
A $W$-loop with base point $x_{0}\in\Sigma$ is defined by    
\begin{enumerate}
\item a sequence $0=t_{0}<\cdots<t_{n}=1,$    
\item continuous paths $c_{i}:[t_{i-1},t_{i}]\rightarrow \Sigma,$ $1\leq i\leq n,$  
\item elements $w_{i}\in W$ defined in a neighborhood of $c_{i}(t_{i})$ for $1\leq i \leq n$ such that $c_{1}(0)=w_{n}c_{n}(1)=x_{0}$ and $w_{i}c_{i}(t_{i})=c_{i+1}(t_{i}),$ where $1\leq i\leq n-1.$
\end{enumerate}

Two $W$-loops are in the same \emph{homotopy class} if one can be obtained from the other by a series of subdivisions, equivalences and deformations.
The homotopy classes of $W$-loops based at $x_{0}\in \Sigma$ form a group $\pi_{1}(W,x_{0})$ called \emph{fundamental group of the pseudogroup} $W$ at the point $x_{0}.$

\begin{rem}
\label{rem-grupo-fundamental-orbifold-Wloopp}
If the orbit space $\Sigma/W$ is a connected orbifold, then $\pi_{1}(W,x)=\pi(\Sigma/W, \rho(x)),$ where $\rho:\Sigma\rightarrow \Sigma/W$ is the natural projection.
\end{rem}

The proof of Theorem \ref{teo-holonomia-trivial} is basically the same proof of Theorem 1.6 of \cite{AlexToeben} apart from the modified Lemma \ref{HolTrivialLema0,5} and Lemma \ref{HolTrivialLema1}. Therefore we only sketch its main steps. For details (e.g. the proofs of Lemma \ref{HolTrivialLema2} and Lemma \ref{HolTrivialLema3}) see \cite{AlexToeben}.

Let $L$ be a regular leaf, $p\in L$ and $\alpha$ a curve in $L$ such that $\alpha(0)=p=\alpha(1).$ 
Let $\zeta(t)$ be the parallel transport of an orthonormal frame $\zeta$ in $p$ along $\alpha$. Note that $\zeta(t)$  is contained in a regular leaf of the singular foliation $\widetilde{\F}$ in $E$. 

We want to show that  $\zeta(0)=\zeta(1).$ 

Since $M$ is simply connected we have a homotopy $G:[0,1]\times [0,1]\to M$ with
\begin{enumerate}
\item $G(0,t)=\alpha(t)$ for all $t\in [0,1]$.
\item $G(s,0)=G(s,1)=p$ for all $s$.
\item $G(1,t)=p$ for all $t$. 
\end{enumerate}
We define $\tilde p:=\zeta(0)$. Let $\pi:E\to M$ be the canonical projection of the transversal frame bundle $E$ of $\F$. We can lift $G$ to a homotopy $\tilde G:[0,1]\times [0,1]\to E$ with
\begin{enumerate}
\item $\tilde G(0,t)=\zeta(t)$ for all $t$.
\item $\tilde G(s,0)=\tilde{p}$ for all $s$.
\item $\tilde G(1,t)=\tilde{p}$ for all $t$.
\item $\pi\circ \tilde G(s,1)=p$ for all $s$.
\end{enumerate}
Let $\Sigma$ be the section of $\F$ that contains $p$ and define $\widetilde\Sigma:=\pi^{-1}(\Sigma)$.
Let $\tilde{\rho}:E\to E/\widetilde\F$ be the natural projection.

\begin{lema}
\label{HolTrivialLema0,5}
Let $\widetilde{L}$ be the lift of a regular leaf $L$. Then $\widetilde{L}$ is closed and embedded in the frame bundle $E_r$.
\end{lema}
\begin{proof}
First we prove that $\widetilde{L}$ is embedded. Let $\tilde{x}\in \widetilde{L}$ and set $x=\pi (\tilde{x})$. Let $U$ be a neighborhood  of $x$ such that $U\cap L$ has only one connected component. Then
\begin{equation}
\label{Eq1-claim-teo-holohomia-trivial-3}
\pi^{-1}(U)\cap\widetilde{L}=\pi^{-1}(U\cap L)\cap\widetilde{L}
\end{equation}

Note that the holonomy of $L$ if finite, because the leaves of $\F$ are closed and embedded. Since the holonomy of $L$ is finite, $\widetilde{L}$ meets $\pi^{-1}(x)$ only a finite number of times. This fact and Equation \ref{Eq1-claim-teo-holohomia-trivial-3} impy that $\widetilde{L}$ is embedded.

Now we have to prove that $\widetilde{L}$ is closed. Let $\{\tilde{x}_{n}\}$ be a sequence contained in $\widetilde{L}$ that converge to a point $\tilde{x}$. Since $L$ is closed we conclude that $\tilde{x}\in\pi^{-1}(L).$
Set $x=\pi(\tilde{x})$. Let $U$ be a neighborhood  of $x$ such that $U\cap L$ has only one connected component. As before, we can note that 
$\pi^{-1}(U)\cap\widetilde{L}$ has only a finite numbers of connected component. This fact and the fact that the sequence $\{\tilde{x}_{n}\}\subset \widetilde{L}$ converge to $\tilde{x}$ imply that  $\tilde{x}\in\widetilde{L}$.

\end{proof}

\begin{lema}
\label{HolTrivialLema1}
 $\tilde{\rho}:\widetilde\Sigma_r\to E_r/\tilde\F$ is a covering map, where  $\widetilde{\Sigma}_r=\widetilde{\Sigma}\cap E_r$.
\end{lema}
\begin{proof}
 
 It is easy to verify that $\tilde{\rho}:\widetilde\Sigma_r\to E_r/\tilde\F$ is surjective.
 Now the result follows from Lemma \ref{HolTrivialLema0,5} and the claim that we will prove below.
 
\begin{claim}
Let $\F$ be a s.r.f.s with trivial holonomy and assume that the leaves are closed and embedded. 
Let $L_{p}$ be a regular leaf and define $\tub_{\epsilon}(L_{p}):=\cup_{x\in L_{p}} D_{x},$ where $D_{x}:=\exp_{x}(B_{\epsilon}(0) )$ for a ball $B_{\epsilon}(0)\subset \nu_{x} L.$ Then there exists an $\epsilon>0$ such that
\begin{enumerate}
\item[a)]for each $x\in L_{p}$ and $y\in D_{x}$ we have $D_{x}\cap L_{y}=\{y\}.$
\item[b)] The map $\exp_{x}:B_{\epsilon}(0)\rightarrow D_{x}$ is a diffeomorphism.
\end{enumerate}
\end{claim}

Since $L$ is embedded and has trivial holonomy we can find an $\epsilon>0$ such that $D_{p}$ is contained in a normal neiborhood of $p$ and  $D_{p}\cap L_{y}=\{y\}$ for each $y\in D_{p}.$

a) Suppose that there are two points $y_1$ and $y_2$ that belong to $D_{x}\cap L_{y}.$ Then there are two vectors $\xi_{1}$ and $\xi_{2}\in\nu L_{x}$ such that $\exp_{x}(\xi_{i})=y_{i}$ for $i=1,2.$ Since the holonomy of $L_{p}$ is trivial, we can extend $\xi_{i}$ to a global parallel normal field. Now using the fact that $\F$ is equifocal and that $D_{p}$ is contained in a normal neiborhood of $p$, we conclude that $\exp_{p}(\xi_{1}(p))$ and $\exp_{p}(\xi_{2}(p))$ are two different points contained in $D_{p}\cap L_{y}$ and this contradicts our choise of $\epsilon.$

b) We conclude that the map $\exp_{x}:B_{\epsilon}(0)\rightarrow D_{x}$ is a bijection, by the same argument used in the proof of Item a).
Therefore it sufficies to prove that the map $\exp_{x}:B_{\epsilon}(0)\rightarrow D_{x}$ is a local diffeomorphism. 
 We can extend  a  vector $\xi\in\nu L_{x}$ to a  global parallel normal field. Gluing germs of holonomy maps along the curve $\gamma_{p}(t)=\exp_{p}(t\xi)$, ($0\leq t\leq 1$). we can construct a local isometry $\varphi:U_{p}\rightarrow U_{x}$ where  $U_x$ (respectively $U_p$) is a neighborhood of the curve $\gamma_x$ (respectively $\gamma_p$) in the section that contain $x$ (repectively $p$).   
The existence of the local isometry $\varphi$ implies that the vector $\xi(x)$ is not a critical point of the map $\exp_{x}^{\perp}$. The arbitrarity of choise of the vector $\xi(x)$ implies that $\exp_{x}^{\perp}|_{B_{\epsilon}(0)}$ is a local diffeomorphism.

\end{proof}

We define $R(s,t)=(1-t,1-s)$ and $\gamma_0:=\tilde G\circ R(0,\cdot)$. $\tilde G\circ R$ is a homotopy between $\gamma_0$ and the constant curve $\gamma_{\tilde p}\equiv \tilde p$. We have $\gamma_0(0)=\tilde p$ and $\pi\circ \gamma_0\equiv p.$ These facts imply  that $\gamma_0$ and $\gamma_{\tilde p}$ are contained in $\widetilde{\Sigma}_{r}.$

It is easy to check that

\begin{lema}
\label{HolTrivialLema1,5}
$\gamma_{0}$ is a $\WSrt$-loop based at $\tilde{p}$ and in particular a $\WSt$-loop based at $\tilde{p}.$
\end{lema}

Using the holonomy map and trivializations of $\widetilde{\F},$ we can project the homotopy $\tilde G\circ R$ to $\WSt$ -loop deformations on $\widetilde{\Sigma}$ and prove the next lemma.

\begin{lema}
\label{HolTrivialLema2}
$\gamma_{0}$ and the trivial $\WSt$-loop $\gamma_{\tilde p}$ belong to the same homotopy class of  $\pi_{1}(\WSt, \tilde{p}).$
\end{lema}

Finally reflecting $W$-loops in the walls of Weyl chambers it is possible to prove the lemma below. 

\begin{lema}
\label{HolTrivialLema3}
Consider two $\WSt$-loops $\delta_{0}$ and $\delta_{1}$ based at $\tilde{p}$ that belong to the same homotopy class of $\pi_{1}(\WSt, \tilde{p}).$ Suppose that $\delta_{0}$ and $\delta_{1}$ are contained in $\widetilde{\Sigma}_{r}.$
Then $\delta_{0}$ and $\delta_{1}$ belong to the same homotopy class of $\pi_{1}(\WSrt, \tilde{p}).$ 
\end{lema}

Lemmas \ref{HolTrivialLema1,5}, \ref{HolTrivialLema2}, \ref{HolTrivialLema3} and the fact that $\pi_{1}(\WSrt,\tilde{p})=\pi_{1}(E_{r}/\widetilde{\F},\rho(\tilde{p}))$ (see Remark \ref{rem-grupo-fundamental-orbifold-Wloopp}) 
imply that $\rho\circ\gamma_0$ and the constant curve $\rho\circ\gamma_{\tilde p}$ are homotopic in $E_r/\tilde\F$ fixing endpoints. The lift of this homotopy along the covering $\tilde{\rho}:\tilde\Sigma_r\to E_r/\F$ (see Lemma \ref{HolTrivialLema1}) to the curve $\gamma_0$ in $\tilde\Sigma_r$ is a homotopy to a constant curve fixing endpoints. Thus $\zeta(0)=\gamma_0(0)=\gamma_0(1)=\zeta(1)$.

\end{proof}
\subsection{Some applications}

\begin{teo}[\cite{AlexToeben}]
\label{teo-holonomia-trivial-secao-flat}
Let $\F$ be a s.r.f.s.\ on a simply connected riemannian manifold $M$. 
Assume that the sections are flat. Then each regular leaf has trivial holonomy.
\end{teo}
\begin{proof}
If the sections of $\F$ are flat, $E$ is a smooth bundle and $\widetilde{\F}$ is a (smooth) singular foliation. Let  $\s:U\rightarrow E$ be the cross-section  with respect to $\tilde q$, which was used in the construction of the bundle $E$ (see  Proposition \ref{prop-fibrado-referenciais}).We can define a distribution $H$ on $E$ by $H_{\tilde{q}}:=T_{\tilde{q}}\s(U).$
It is not difficult to check that this distribution is integrable. This implies that $E$ is foliated by submanifolds $\{\widetilde{M}_{\tilde{x}}\}$ and for each $\tilde x\in E$ the map $\pi:\widetilde{M}_{\tilde{x}}\rightarrow M$ is a covering map. For each manifold $\widetilde{M}_{\tilde{x}}$ the lift of $\F$ along $\pi$ coincides with $\widetilde{\F}|\widetilde{M}_{\tilde{x}}.$ This is exactly what happens in Example \ref{trivialex}.
The covering map $\pi:\widetilde{M}_{\tilde{x}}\rightarrow M$ is a diffeomorphism if $M$ is simply connected. This implies that the regular leaves of $\F$ have trivial holonomy.
\end{proof}

Note that, in the above result, we do not assume that the leaves are embedded or closed.

\begin{cor}
\label{cor-holonomia-grupo-fundamental}
Let $\F$ be a s.r.f.s on a complete riemannian manifold $M.$ Suppose that the cardinality of the fundamental group $\pi_{1}(M)$ is equal to $n.$ Assume that one of the two conditions below is satisfied
\begin{enumerate}
\item The leaves of $\F$ are embedded and closed.
\item The sections are flat.
\end{enumerate}
Then the cardinality of the holonomy of $\F$ is lower or equal to $n.$ 
\end{cor}

\begin{proof}
Let $\widetilde{M}$ be the riemannian covering space of $M$ and $\pi:\widetilde{M}\rightarrow M$ be the riemannian covering map. Denote $\widetilde{\F}$ as the lift of the foliation $\F$. Let $x_{0}$ be a regular point and consider a loop $\beta\subset L_{x_0}$ with $\beta(0)=x_{0}=\beta(1)$. Finally define $\tilde{\beta}$ as the leaft of $\beta$ such that $\tilde{\beta}(0)=x_{0}.$ We claim that
\begin{equation}
\label{Eq-levantamento-holomomia-1}
\varphi_{[\beta]}\circ \pi=\pi\circ\tilde{\varphi}_{[\tilde{\beta}]}
\end{equation}

In fact we can find a partition  $0=t_{0}<\cdots<t_{n}=1$ such that $\beta_{i}:=\beta|_{[t_{i-1,i},t_{i}]}$ is contained in 
a distinguished neighorhood of a foliation chart of $\F$. We can also assume that $\beta_{i}$ is contained in a neighborhood $U$ such that $\pi^{-1}(U)$ is a disjoint union of open subsets $U_{\alpha}$ such that $\pi: U_{\alpha}\rightarrow U$ is a diffeomorphism. Clearly
\begin{equation}
\label{Eq-levantamento-holomomia-2}
\varphi_{[\beta_{i}]}\circ \pi=\pi\circ\tilde{\varphi}_{[\tilde{\beta}_{i}]}.
\end{equation}
for each $i.$ 
Assume by induction that 
\begin{equation}
\label{Eq-levantamento-holomomia-2-5}
\varphi_{[\beta_{i}\circ\cdots\circ\beta_{1}]}\circ \pi=\pi\circ\tilde{\varphi}_{[\tilde{\beta}_{i}\circ\cdots\circ\tilde{\beta}_{1}]}.
\end{equation}
Therefore 
\begin{eqnarray*}
\varphi_{[\beta_{i+1}\circ\cdots\circ\beta_{1}]}\circ \pi& = & \varphi_{[\beta_{i+1}]}\circ \varphi_{[\beta_{i}\circ\cdots\circ\beta_{1}]}\circ \pi\\
&\stackrel{\mathrm{Eq.}\ref{Eq-levantamento-holomomia-2-5}} {=} & \varphi_{[\beta_{i+1}]}\circ \pi\circ\tilde{\varphi}_{[\tilde{\beta}_{i}\circ\cdots\circ\tilde{\beta}_{1}]}\\
 &\stackrel{\mathrm{Eq.}\ref{Eq-levantamento-holomomia-2}} {=}&  \pi\circ\tilde{\varphi}_{[\tilde{\beta}_{i+1}]} \circ \tilde{\varphi}_{[\tilde{\beta}_{i}\circ\cdots\circ\tilde{\beta}_{1}]}\\ 
&=& \pi\circ \tilde{\varphi}_{[\tilde{\beta}_{i+1}\circ\cdots\circ\tilde{\beta}_{1}]},
\end{eqnarray*}
and this prove  Equation \ref{Eq-levantamento-holomomia-1}.

On the other hand, it follows from  Theorem \ref{teo-holonomia-trivial} and Theorem \ref{teo-holonomia-trivial-secao-flat} that the holonomy of $\widetilde{\F}$ is trivial. Thus there exist only $n-1$ holonomy $\varphi_{[\tilde{\beta}]}$ between $\tilde{x}_{0}$ and the others points $\tilde{x}_{1}\ldots\tilde{x}_{n-1} \in \pi^{-1}(x_{0})$. This fact and Equation \ref{Eq-levantamento-holomomia-1} imply the result.

\end{proof}

In \cite{AlexToeben} T\"{o}ben and I proved the existence of fundamental domains in each section of a s.r.f.s. when the leaves are compact and $M$ is simpy connected. Due to Theorem \ref{teo-holonomia-trivial} we can reformulate our result as follows. 


\begin{teo}
\label{teo-fundamental-domain}
Let $\F$ be a s.r.f.s.\ on a simply connected riemannian manifold $M$. 
Suppose also that the leaves of $\F$ are  closed embedded. Then
\begin{enumerate}
\item $M/ \F$ is a simply connected Coxeter orbifold.
\item Let $\Sigma$ be a section of $\F$ and 
let $\Pi:M\rightarrow M/\F$ be the canonical projection. 
Denote by $\Omega$ a connected component of the set of regular points 
in 
$\Sigma$. Then  $\Pi:\Omega\rightarrow M_{r}/\F$ and 
$\Pi:{\overline{\Omega}}\to M/\F$ are homeomorphisms, where $M_r$ 
denotes the set of regular points in $M$. In addition, 
$\Omega$ is convex, i.e.~for any two points~$p$ and~$q$ in~$\Omega$, 
every minimal geodesic segment between $p$ and $q$ lies 
entirely in $\Omega$.
\end{enumerate}
\end{teo}

The existence of fundamental domain turns out to be a useful tool to study s.r.f.s. Indeed, this was one of the techniques used by  Gorodski and I in \cite{AlexGorodski}  to prove that the leaves of a s.r.f.s are pre image of a transnormal map, when the leaves are compact, the sections are flat and $M$ is simply connected (for the definition of transormal map see Definition \ref{dfn-transnormal-map}).  
Due to Theorem \ref{teo-fundamental-domain}, we can reformulate our result as follows.

\begin{teo}
\label{teo-s.r.f.s-transnormal}
Let $\mathcal F$ be a singular riemannian foliation with sections
on a complete simply connected riemannian manifold $M$.
Assume that the leaves of $\mathcal F$ are closed embedded and that  
$\mathcal F$ admits a flat section of dimension $n$. 
Then the leaves of $\mathcal F$ are given by 
the level sets of a transnormal map $F:M\to\R^n$.
\end{teo}

The above theorem generalizes 
previous results of Carter and West~\cite{CarterWest2},
Terng~\cite{terng} and Heintze, Liu and Olmos~\cite{HOL}
for isoparametric submanifolds. 
It can also be viewed as a converse to the main result 
in~\cite{Alex1}, and as a global version of Proposition \ref{prop-placas-transnormais}


\section{Blow-up and extension of the holonomy map}

In \cite{Toeben} T\"{o}ben  used the blow up technique to study  equifocal submanifolds (which he called submanifold with parallel focal structure). He gave a necessary and sufficient condition for a closed embedded  equifocal submanifold to induce a s.r.f.s (see \cite{Alex3} for an alternative proof). 

The aim of this section is twofold. First we  we will briefly recall some of T\"{o}ben's results about blow-up of s.r.f.s (see Theorem \ref{Th-Toeben} and Proposition \ref{Prop-Toeben}). Then
we will use this technique to get conditions under which a holonomy map can be extended to a global isometry (see Proposition \ref{Prop-extensao}). 

We start by recalling the blow-up technique.

\begin{teo}[T\"{o}ben \cite{Toeben}]
\label{Th-Toeben}
Let $\F$ be a s.r.f on a complete riemannian manifold $M.$ Then 
\begin{enumerate}
\item[a)] Set $\widehat{M}:=\{T_{p}\Sigma |\  p\in N, \Sigma$ is a section of $\F$ through $p\}.$ Then $\widehat{M}$ carries  a natural differentiable structure, for which the inclusion into the Grassmann bundle $G_{k}(TM)$ is an immersion. Moreover, $\widehat{M}$ has a natural 
 riemannian/ totally geodesic bifoliation $(\widehat{\F}, \widehat{\F}^{\perp})$, with respect to the pull-back metric.  We have $\widehat{\F}^{\perp}=\{T\Sigma| \ \Sigma$ is a section of $\F\}$.
\item[b)] The footpoint map $\hat{\pi}: (\widehat{M},\widehat{\F})\rightarrow (M, \F)$ is foliated and maps each horizontal leaf of $\widehat{F}^{\perp}$ isometrically to the corresponding section $\Sigma$ of $\F.$
\end{enumerate}

\end{teo}

\begin{rem}
The result above is a strengthening of Boualem's result \cite{Boualem}. He stated it for some differentiable structure and some metric. T\"{o}ben proved it for the natural differential structure and natural metric. Morover he did not need that the leaves were relatively compact, as assumed by Boualem. 
\end{rem}

In order to study the singular holonomy of $\F$, T\"{o}ben considered the universal covering space of $\widehat{M}$, which turns out to be diffeomorphic to $\tilde{L}\times\tilde{\Sigma}$ as we recall below.

\begin{lema}[\cite{BH1}]
\label{Lemma-Toeben}Let $(\F,\F^{\perp})$ be the bifoliation on $\widehat{M}$ defined in Theorem \ref{Th-Toeben}.
Let $\hat{x}_0$ a point of $\widehat{M}$   , $\hat{\beta}:[0,1]\rightarrow \widehat{M}$ be a curve contained in the leaf  $\widehat{L}_{\hat{x}_{0}}\in\widehat{\F}$ and $\hat{\gamma}:[0,1]\rightarrow  
\widehat{M}$ a curve contained in the leaf  
$\widehat{\Sigma}_{\hat{x}_{0}}\in \widehat{F}^{\perp}$ such that $\hat{\gamma}(0)=\hat{\beta}(0)=\hat{x}_{0}$. Then there exists a unique continuou map $\widehat{H}=\widehat{H}_{(\hat{\beta},\hat{\gamma})}:[0,1]\times[0,1]\rightarrow \widehat{M}$ with
\begin{enumerate}
\item $\widehat{H}(\cdot,0)=\hat{\beta}.$
\item $\widehat{H}(0,\cdot)=\hat{\gamma}.$
\item $\widehat{H}(\cdot, t)$ is contained in a leaf of $\widehat{\F}.$
\item $\widehat{H}(s,\cdot)$ is contained in a leaf of $\widehat{\F}^{\perp}.$
\end{enumerate}
\end{lema}
The continuous map $\widehat{H}$ is called \emph{rectangle} with \emph{initial vertical} (respectively \emph{horizontal}) curve $\hat{\beta}$ (respectively $\hat{\gamma}$). 
\begin{rem}
\label{rem-Lemma-Toeben}
For a curve $\beta:[0,1]\rightarrow M$ in a regular leaf of $\F$ and a curver $\gamma:[0,1]\rightarrow M$ in a section, both starting  in a regular point $x_{0}$, we can define the lift $\hat{\beta}(t):=T_{\beta(t)}\Sigma_{\beta(t)}$ and $\hat{\gamma}(t):=T_{\gamma(t)}\Sigma_{\gamma(0)}$.
Clearly $\hat{\pi}\circ\hat{\beta}=\beta$ and $\hat{\pi}\circ\hat{\gamma}=\gamma.$
As remarked in \cite{Toeben}, the above lemma is also true for a s.r.f.s $\F$. If we write $H_{(\beta,\gamma)}$ for the rectangle with  initial vertical (respectively horizontal) curve $\beta$ (respectively $\gamma$) we can note that $H_{(\beta,\gamma)}=\hat{\pi}\circ \widehat{H}_{(\hat{\beta},\hat{\gamma})}$. 

\end{rem}

We recall that the universal cover $\widetilde{M}$ of a manifold $M$ is equal to the set of equivalence of curves starting from a fixed point $x_{0}$, where the equivalence is given by homotopy fixing endpoints. Therefore, for a regular point $x_0$ we have 
\[\widetilde{L}=\{[\beta] \  |  \ \beta \mathrm{ \ is \ vertical \ and \ } \beta(0)=x_0\} ,\]
\[\widetilde{\Sigma}=\{[\gamma] \ | \ \gamma \mathrm{\ is \ horizontal \ and \ } \gamma(0)=x_0\},\]
\[\widetilde{\widehat{M}}=\{[\mu] \ | \ \mu \mathrm{\ is \ a \ curve \ in \ \widehat{M} \ and \ } \mu(0)=\hat{x}_{0}=T_{x_0}\Sigma_{x_{0}}\}. \]

Now consider the manifold $\widetilde{L}\times \widetilde{\Sigma}$ provided with the natural bifoliation and the covering space $\widetilde{\widehat{M}}$ provided with the pull-back bifoliation of the covering map $\widetilde{\widehat{M}}\rightarrow \widehat{M}.$ It follows from Blumenthal and Hebda \cite{BH2} that the map $\Phi:\widetilde{L}\times\widetilde{\Sigma}\rightarrow \widetilde{\widehat{M}},$ defined as 
$([\beta],[\gamma])\rightarrow [t\rightarrow \widehat{H}_{(\hat{\beta},\hat{\gamma})}(t,t)] ,$
is a bifoliated diffeomorphism (i.e., foliated with respect to both pairs of foliations). We conclude that
\[
\begin{array}{llcl}
\Psi:& \widetilde{L}\times\widetilde{\Sigma}&\rightarrow& \widehat{M}\\
     &  ([\beta],[\gamma])    &\rightarrow&  \widehat{H}_{(\hat{\beta},\hat{\gamma})}(1,1)
\end{array}\]
is a bifoliated universal covering map of $\widehat{M}$. Define $\psi: \widetilde{L}\times\widetilde{\Sigma}\rightarrow M$ as $\psi:=\hat{\pi}\circ\Psi.$

The above discution  and Theorem \ref{Th-Toeben} imply the next proposition.

\begin{prop}[\cite{Toeben}]
\label{Prop-Toeben}
The map $\Psi$ is the universal covering map, and it is bifoliated with respect to the natural bifoliation of $\widetilde{L}\times\widetilde{\Sigma}$ and to $(\widehat{M};\widehat{F},\widehat{F}^{\perp}).$ The map $\psi$ is foliated with respect to the vertical foliation on $\widetilde{L}\times\widetilde{\Sigma}$ and $(M;\F)$, and its restriction to a horizontal leaf is a riemannian covering to a section. 
\end{prop}

\begin{cor}[\cite{Toeben}]
Let $\F$ be a s.r.f.s on a complete riemannian manifold $M.$ Then the sections have the same riemannian universal cover. Similarly the regular leaves of $\F$ have the same universal cover. 
\end{cor}

Using Theorem \ref{Th-Toeben} and Proposition \ref{Prop-Toeben},  T\"{o}ben proved the next result.

\begin{prop}[\cite{Toeben}]
Let $\F$ be a s.r.f.s on a complete riemannian manifold $M.$ Assume that the sections are embedded. Then there exists a section $\Sigma$ such that its Weyl pseudogroup $W_{\Sigma}$ is in fact a group.
\end{prop}

Now we  infer from  Theorem \ref{Th-Toeben} and Lemma \ref{Lemma-Toeben} that the each holonomy map can be extended to a global isometry, if the sections are embedded and if there exists a singular point $q$ such that the leaf passing through this point is just $q.$  

\begin{prop}
\label{Prop-extensao}
Let $\F$ be a s.r.f.s on a complete riemannian manifold $M.$ Assume that the sections of $\F$ are embedded and that there is a leaf which is point $q$, i.e., $L_{q}=\{q\}$. Let $\varphi_{[\beta]}:\sigma_{0}\rightarrow\sigma_{1}$ be a holonomy map where $\sigma_{0}$ (respectively $\sigma_{1}$) is a local section of a section $\Sigma_{0}$ (respectively $\Sigma_{1}$). Then there exists an isometry $\varphi:\Sigma_{0}\rightarrow\Sigma_{1}$ such that $\varphi|_{\sigma_{0}}=\varphi_{[\beta]}$ and $\varphi(x)=L_{x}.$ 
\end{prop}

\begin{proof}
First we note that the leaves of $\widehat{F}^{\perp}$ are embedded, since the sections of $\F$ are embedded. Then we define $\widehat{L}_{\hat{q}}:=(\hat{\pi})^{-1}(q).$  The fact that $L_{q}=\{q\}$ implies that each leaf of $\widehat{\F}^{\perp}$ meets $\widehat{L}_{\hat{q}}$ once and only once. This implies that the holonomy of the leaves of $\widehat{\F}^{\perp}$ are trivial. Now the result follows from the lemma below.
\begin{lema}
If the leaves of $\widehat{\F}^{\perp}$ are embedded and have trivial holonomy, then each holonomy map of $\F$ admits an extension to a global isometry.
\end{lema}
\begin{proof}
It sufficies to prove that for each loop   $\gamma\subset\Sigma_{x_{0}}$ with $\gamma(0)=\gamma(1)=x_{0}:=\beta(0)$, the continuation of the germ of $\varphi_{[\beta]}$ along $\gamma$ leads back to the initial germ.

Let $\hat{\beta}$ (respectively $\hat{\gamma}$)  be the lift of the curve $\beta$ (respectively $\gamma$) defined in Remark \ref{rem-Lemma-Toeben}. 
Note that $\hat{\gamma}$ is also a loop, because there exists only one point $\hat{x}_{0}$ such that $\hat{\pi}(\hat{x}_{0})=x_{0}$. 
Let $\widehat{H}$ denote the rectangle with initial vertical (respectively horizontal) curve $\hat{\beta}$ (respectively $\hat{\gamma}$) (see Lemma \ref{Lemma-Toeben}).

We can find partitions $0=s_{0}<\cdots < s_{n}=1$ and $0=t_{0}<\cdots t_{n}=1$ so that $\widehat{H}|_{[s_{i-1},s_{i}]\times[t_{j-1},t_{j}]}$ is contained in a distinguished neighorhood of a foliation chart of $\widehat{\F}.$ Set $\beta_{i}=\beta|_{[s_{i-1},s_{i}]}$ and let $\Sigma_{i}$ denote the section which contains $\beta(s_{i})$.

Define a holonomy map $\varphi^{j}_{1}:\sigma_{0\,j}\rightarrow \sigma_{1\, j}$ such that 
\begin{enumerate}
\item $\varphi^{0}_{1}=\varphi_{[\beta_{1}]},$  
\item $\varphi^{j}_{1}|_{\sigma_{0\, j-1}\cap\sigma_{0\,j}}=\varphi^{j-1}_{1}|_{\sigma_{0\, j-1}\cap\sigma_{0\,j}},$ 
\end{enumerate}
where $\sigma_{0\,j}$ are local sections of $\Sigma_{0}$ centered at $\gamma(t_{j})$ and $\sigma_{1\,j}$ are local sections of $\Sigma_{1}.$

We want to prove 
\begin{equation}
\label{extensao-primeira}
\varphi^{n}_{1}|_{\sigma_{0\, n}\cap\sigma_{0\,0}}=\varphi^{0}_{1}|_{\sigma_{0\, n}\cap\sigma_{0\,0}}.
\end{equation}

We note that each holonomy map $\varphi^{j}_{1}$ is  associated to an holonomy map $\hat{\varphi}^{j}_{1}:\hat{\sigma}_{0\,j}\rightarrow\hat{\sigma}_{1\,j}$ of the regular foliation $\widehat{\F},$ where $\hat{\sigma}_{i\,j}$ is a neigborhood of $\widehat{\Sigma}_{i}$ (the leaf of $\widehat{\F}^{\perp}$ which contains $\hat{\beta}(s_{i})$) such that

\begin{enumerate}
\item $\hat{\pi}(\hat{\sigma}_{i\,j})=\sigma_{i\,j}$,
\item $\hat{\pi}\circ\hat{\varphi}^{j}_{1}=\varphi^{j}_{1}\circ\hat{\pi}$.
\end{enumerate}
Note that $\hat{\sigma}_{0\, n}\cap\hat{\sigma}_{0\,0}\neq\emptyset$ since $\hat{\gamma}$ is a loop.

To prove Equation \ref{extensao-primeira} it sufficies to prove 
\begin{equation}
\label{extensao-segunda}
\hat{\varphi}^{n}_{1}|_{\hat{\sigma}_{0\, n}\cap\hat{\sigma}_{0\,0}}=\hat{\varphi}^{0}_{1}|_{\hat{\sigma}_{0\, n}\cap\hat{\sigma}_{0\,0}}.
\end{equation}

Now Equation \ref{extensao-segunda} follows direct from te fact that the holonomy of $\widehat{\F}^{\perp}$ is trivial and that
$\widehat{H}|_{[s_{i-1},s_{i}]\times[t_{j-1},t_{j}]}$ is contained in a distinguished neighorhood of a foliation chart of $\widehat{\F}.$

By induction we can prove that 
\begin{equation}
\label{extensao-terceira}
\varphi^{n}_{i}|_{\sigma_{0\, n}\cap\sigma_{0\,0}}=\varphi^{0}_{i}|_{\sigma_{0\, n}\cap\sigma_{0\,0}}.
\end{equation}
Defining $\varphi^{j}=\varphi^{j}_{n}\circ\cdots\circ\varphi^{j}_{1}$, we conclude that  
\begin{equation}
\label{extensao-quarta}
\varphi^{n}|_{\sigma_{0\, n}\cap\sigma_{0\,0}}=\varphi^{0}|_{\sigma_{0\, n}\cap\sigma_{0\,0}}= \varphi_{[\beta]}|_{\sigma_{0\, n}\cap\sigma_{0\,0}}.
\end{equation}

Equation \ref{extensao-quarta} implies that the continuation of the germ of $\varphi_{[\beta]}$ along $\gamma$ leads back to the initial germ. This completes the proof. 
\end{proof}

\end{proof}

\bibliographystyle{amsplain}

\end{document}